\documentclass[11pt]{article}
\usepackage{graphicx,amsmath,bm, amsthm,mathrsfs,amssymb, braket, verbatim}
\usepackage{caption}
\usepackage{subcaption}
\usepackage[usenames]{color}
\usepackage{ulem,mathtools}
\usepackage{pdfpages}
\usepackage{lscape}
\usepackage{cite}
\usepackage{authblk}
\usepackage{cite}
\usepackage{tcolorbox}
\usepackage[section]{placeins}
\usepackage{placeins}

\newcommand{\ee}{\end{equation}}

\newcommand{\reff}[1]{(\ref{#1})}
\newcommand{\beq}{\begin{equation}}
\newcommand{\eeq}[1]{\label{#1}\end{equation}}
\newcommand{\beqa}{\begin{eqnarray}}
\newcommand{\eea}{\end{eqnarray}}
\newcommand{\eeqa}[1]{\label{#1}\end{eqnarray}}
\newcommand{\beg}{\begin{equation*}}
\newcommand{\eeg}{\end{equation*}}

\newcommand{\eq}{\!=\!}
\newcommand{\p}{\!+\!}
\newcommand{\m}{\!-\!}

\newcommand{\bsplit}{\begin{split}}
\newcommand{\esplit}{\end{split}}

\usepackage{circuitikz} 
\usepackage[capposition=bottom]{floatrow} 
\usepackage{epigraph} 
\usepackage{appendix}
\usepackage{rotating}
\allowdisplaybreaks

\title{Summing the reciprocal of the polynomial \\appearing in Fermat's Last Theorem}
\author[]{Ariel Edery\thanks{aedery@ubishops.ca}}
\affil[]{Department of Physics and Astronomy, Bishop's University, 2600 College Street, Sherbrooke, Qu\'{e}bec, Canada, J1M 1Z7.\vspace{1em}}
\begin{document}
\date{}
\maketitle
\begin{abstract}
Consider the polynomial $f \eq x^N \p y^N \m z^N$ where $x,\,y$ and $z$ are positive integers and $N \ge 3$ is an integer. By Fermat's Last Theorem, $f$ is never zero so that its reciprocal, $1/f$, encounters no singularities. We therefore study the finite sum of the reciprocal:\\ $S(m,N)\!=\!\sum_{x=1}^m\sum_{y=1}^m\sum_{z=1}^{m} \frac{1}{f}$. The terms $1/f$ can be positive, negative and their magnitude are bounded by unity i.e. $|1/f|\le 1$. A key observation is that $S(m,N)$ can be split into two convenient parts: a dominant contribution $D(m,N)$ that has a simple analytical expression and a remainder $R(m,N)$ which is more complicated but negligible compared to $D(m,N)$. Therefore, $S(m,N)$ is almost identical to $D(m,N)$. The analytical expression for 
$D(m,N)$ is $(2\, m -1) \,H_m^{(N)}$ where $H_m^{(N)}=\sum_{x=1}^m \frac{1}{x^N}$. After summing just a few terms, $H_m^{(N)}$ approaches quickly the Riemann zeta function $\zeta(N)$. Therefore, the original sum $S(m,N)$ has a remarkably simple expression: it is basically linear in $m$ with slope equal to $2\,\zeta(N)$. Its linear behavior is not an asymptotic result; plots of $S(m,N)$ vs. $m$ for different $N$ show a straight line starting at $m=1$. $S(m,N)$ deviates slightly from a straight line over a small interval $8\le m \le 12$ for the case $N=3$. This slight deviation can be traced to Fermat near misses where $x^3+y^3-z^3=\pm 1$ (for $z\ne x$ and $z\ne y$); these create a jump in the remainder $R(m,3)$ at $m=9$. We make a numerical and analytical study of the remainder $R(m,N)$. From our numerical analysis, one could safely conclude that $R(m,N)$ converges for $N\ge 4$ but it was harder to tell whether the case $N=3$ converged or diverged. An analytical study based on a comparison of $R(m,N)$ to its Cauchy principal value integral, shows that $R(m,3)$ likely diverges logarithmically. The analytical study also showed that $R(m,N)$ converges for $N\ge 4$ in agreement with our numerical analysis. We discuss in the conclusion some interesting questions for future investigation.     
\end{abstract}
\setcounter{page}{1}
\newpage
\section{Introduction}\label{Intro}
The polynomial function $f=x^N + y^N - z^N$ where $x$, $y$ and $z$ are positive integers and $N\ge 3$ is an integer, has the special property of never being equal to zero. This is due to Fermat's Last Theorem which was proved in 1995 \cite{Wiles, Wiles2}. This implies that its reciprocal $\frac{1}{f}$ does not encounter a singularity. It can be positive or negative and its magnitude is less than or equal to unity since $|f|\ge 1$. In this work we consider its sum over $x$, $y$ and $z$ i.e.
\beq
S(m,N)=\sum_{x=1}^m\sum_{y=1}^m\sum_{z=1}^m \frac{1}{x^N+y^N-z^N}
\eeq{Sum} 
where $m>1$ is an integer. The sum does not have the simplifying feature of having a summand which is positive definite. The negative coefficient in front of $z^N$ allows for negative contributions besides the positive ones. Moreover, the polynomial in the denominator is of a general degree $N$ and one might expect very different behaviors depending on the value of the integer $N$. All of this would make it seem like the sum $S(m,N)$ would yield a complicated result and one where different values of $N$ might even have to be considered separately. Yet, remarkably, this sum is linear in $m$ and can be approximated to very high accuracy by the simple analytical result $(2\,m-1) \,\zeta(N)$. In short, a plot of $S(m,N)$ vs. $m$ yields a straight line with slope $2\, \zeta(N)$. Moreover, this is \textit{not} an asymptotic result but valid throughout; the straight line starts at $m=1$. How does a simple analytical result stem from such a complicated sum? It turns out that the sum $S(m,N)$ can be separated into two parts: a dominant contribution $D(m,N)$ and a remainder $R(m,N)$. $D(m,N)$ can be evaluated analytically and this is where the linearity in $m$ stems from. The remainder $R(m,N)$ is a complicated sum but it turns out that its value is negligible compared to the linear term i.e. the ratio $|R(m,N)/D(m,N)|$ is a lot less than unity and approaches zero as 
$m$ increases. So in essence, the complicated part of the sum is relegated to a negligible remainder. 

A numerical plot of $R(m,N)$ vs. $m$ for $N=3,4,5$ and $6$ shows that the case $N=3$ fluctuates the most and that $N=4$ fluctuates considerably also. The cases $N=5$ and $N=6$ had almost no fluctuations and converged quickly. Our numerical analysis could not confirm whether $R(m,3)$ diverges or not. It showed that if it diverged, the growth would be extremely slow, non-monotonic with considerable local fluctuations. If it converged, this would only appear at values of $m$ that are orders of magnitude beyond those we computed. An analytical study comparing the sum $R(m,N)$ to the Cauchy principal value of its integral showed that $R(m,3)$ likely diverges logarithmically. This does not alter the fact that the ratio $|R(m,3)/D(m,3)|$ is a lot less than unity and tends to zero as $m$ increases. $S(m,3)$ is a straight line like all the $N>3$ cases. It was observed that $N=4$ had fluctuations over a significant range but ultimately converged at higher $m$ values. Our analytical study showed that $R(m,N)$ converges for $N\ge 4$ in agreement with our numerical analysis. There is a jump in the remainder for $N=3$ of significant magnitude near $m=9$. This is why $S(m,3)$ deviates from a straight line over a small interval $8\le m\le 12$. We explain why this occurs at $N=3$ but not for any higher $N$ values.       
 
There is substantial literature on the sum of a reciprocal of a polynomial involving integers to various powers. Here we quote a few articles on this topic that are most relevant to our work. Probably the best known of these reciprocal sums is the Epstein zeta function $Z_Q(s)=\sum_{(m,n)\ne(0,0)}\,1/Q(m,n)^s$ where $Q(m,n)$ is a positive definite quadratic form and $\Re(s)>1$. A thorough detailed study of the values of Epstein zeta functions for rational integers $s=k>2$ was carried out in \cite{Smart} where $Q(m,n)$ was given by $a\,m^2+2\,b\,m\,n +c\,n^2$ with $a>0$ and $a\,c-b^2=1$. In \cite{Tsumura} the author considers multiple harmonic series of Mordell–Tornheim type, explicitly described as partial sums of Mordell–Tornheim zeta series and proves some reduction results. There is also recent work involving mixed zeta functions \cite{Essa}. The Epstein-zeta function and generalizations thereof have also been studied in the physics context of Casimir energies where they provide regularization of an otherwise infinite sum \cite{Elizalde1,Elizalde2,Cheng,Edery}.   

In our work, the denominator can be positive or negative because it avoids being zero due to Fermat's Last Theorem. The polynomial is therefore not required to be positive definite for its reciprocal to be non-singular. Most of the sums of the reciprocal of a polynomial that have been studied are infinite sums of a positive definite quantity. There is a good reason for this: the series converges. In our case, we have a finite 
(triple) sum where the three integer variables run from $1$ up to $m$ inclusively. Though the sum $S(m,N)$ diverges as $m\to \infty$, the interesting thing to study here is its dependence on the finite integers $m$ and $N$. One might think that the negative sign in front of $z^N$ would make things complicated but in fact it leads to a simple analytical result. Though $N$ has a significant effect on the remainder, it plays a minor role in the linear term: it appears in the slope $2\,\zeta(N)$ which does not have a very strong dependence on $N$; the slope has a maximum value of $2.40$ at $N=3$ and decreases rapidly to $2$ as $N$ increases. 

Our paper is organized as follows. In section 2 we split $S(m,N)$ into a dominant contribution $D(m,N)$ and a remainder $R(m,N)$. We show $D(m,N)$ is linear in $m$ with slope $2\,\zeta(N)$. We plot $S(m,N)$ vs. $m$ for $N=3,4$ and $5$ and show that one obtains a straight line in all cases starting at small $m$. In section 3 we analyze the remainder term $R(m,N)$. We study it numerically in section 3.1 and provide plots and table of values for $N=3,4,5$ and $6$. In section 3.2 we make an analytical study of the convergence or divergence of $R(m,N)$ by comparing it to the Cauchy principal value of its integral. Section 4 is the conclusion where we summarize our main results and discuss questions worth investigating in the future.  
 
\section{Summing the inverse of $x^N +y^N-z^N$: linear term plus remainder} 

Consider the polynomial function $f=x^N +y^N-z^N$ where $x, y$ and $z$ are positive integers and $N$ is an integer greater than $2$. By Fermat's Last Theorem, $f$ is never zero and therefore its reciprocal never hits a singularity. The reciprocal is finite and ranges between $1$ and $-1$ inclusively since $|f|\ge 1$. We therefore define $S(m,N)$ as the finite (triple) sum over $x$, $y$ and $z$ of the inverse of $f$: 
\beq
S(m,N) =\sum_{x=1}^m\,\sum_{y=1}^m \,\sum_{z=1}^m \dfrac{1}{x^N+y^N-z^N}=
\sum_{x,y,z=1}^m \dfrac{1}{x^N+y^N-z^N}
\eeq{Smn}
where $m$ is an integer greater than $1$ and the triple sum has been expressed using one sum symbol for notational simplicity. The key observation is that the dominant contribution to this sum comes from the cases where $z=x$ and $z=y$. For these cases, the triple sum reduces to a double sum but over a single variable. This implies the double sum includes a sum over unity that contributes a factor of $m$. For example, for $z=x$, one is left with $y^N$ in the denominator and hence with the single variable $y$. The double sum includes a sum over $y$ and $x=z$ but the sum over $x=z$ is a sum over unity and simply yields a factor of $m$. The same reasoning applies to the case $z=y$. These two cases lead essentially to a linear dependence on $m$. It is therefore convenient to split the original sum \reff{Smn} into two parts:
\begin{align}
S(m,N) = D(m,n) +R(m,n)
\label{Split}
\end{align}
where $D(m,n)$ is the dominant contribution defined as 
\beq
D(m,N)=\sum_{\substack{x,y,z=1\\z=x\,;\,z=y}}^m \dfrac{1}{x^N+y^N-z^N}
\eeq{Dmn}
and $R(m,N)$ is the remainder defined as 
\beq
R(m,N)=\sum_{\substack{x,y,z=1\\z\ne x\,;\,z\ne y}}^m \dfrac{1}{x^N+y^N-z^N}\,.
\eeq{Rmn}   
We now evaluate $D(m,N)$. We first perform the sum when $z=x$ and then the sum when $z=y$. The sum when $x=y=z$ is therefore counted twice and we therefore subtract one such sum. This yields     
\begin{align}
D(m,N)&=\sum_{x=1}^m\sum_{y=1}^m \dfrac{1}{y^N} + \sum_{y=1}^m\sum_{x=1}^m \dfrac{1}{x^N}
- \sum_{x=1}^m \dfrac{1}{x^N}\nonumber\\
&=m\,\sum_{y=1}^m \dfrac{1}{y^N} + m\,\sum_{x=1}^m \dfrac{1}{x^N}
- \sum_{x=1}^m \dfrac{1}{x^N}\nonumber\\
&= (2\,m-1)\,\sum_{x=1}^m \dfrac{1}{x^N} \nonumber\\
&= (2\,m-1)\,H_m^{(N)}  
\label{Linear} 
\end{align}
where $H_m^{(N)}$ is defined as
\beq
H_m^{(N)} =\sum_{x=1}^m \dfrac{1}{x^N}\,.
\eeq{Hmn}
$H_m^{(N)}$ is finite for any $N>1$ and is called a generalized harmonic number. At 
$m=1$ it is equal to unity and as $m$ increases it approaches rapidly the Riemann zeta function $\zeta(N)$. More specifically, it has the series expansion 
\beq
H_m^{(N)}= \zeta(N) -\frac{1}{(N-1)\,m^{N-1}}+\frac{1}{2\,m^{N}}+\mathcal{O}(1/m^{N+1})\,.
\eeq{SeriesHmn}
Therefore $D(m,N)=(2\,m-1)\,H_m^{(N)}$ is essentially linear in $m$ with a slope equal to $2\,\zeta(N)$. This is not a result limited to the asymptotic regime. $D(m,N)$ is a straight line throughout, right from the start, not just asymptotically.  

Even when we replace $H_m^{(N)}$ by $\zeta(N)$ so that we represent $D(m,N)$ by $(2\,m-1)\,\zeta(N)$, it reproduces almost exactly the full sum $S(m,N)$ given by \reff{Smn}. The remainder $R(m,N)$ is therefore negligible compared to $D(m,N)$. To illustrate this, we plot on the same graph the full sum $S(m,N)$ and the quantity $(2\,m-1)\,\zeta(N)$ to represent the straight line $D(m,N)$. We generate plots for $N=3$, $N=4$ and $N=5$ in figures 1, 2 and 3 respectively. The two functions overlap so completely that it is hard to distinguish them on the plot. For $N=3$, the full sum $S(m,3)$ follows a straight line everywhere but deviates slightly from it in a small interval $8\le m \le 12$. This is not due to replacing $H_m^{(N)}$ by $\zeta(N)$. We will see in the next section that the deviation is caused by a jump in the remainder $R(m,3)$ over the same interval due to so-called Fermat near misses. For $N=4$ and $N=5$, $S(m,N)$ does not deviate from a straight line at all, as can be seen in the plots of figures 2 and 3 respectively.   
\begin{figure}[t]
	\centering
		\includegraphics[scale=1.2]{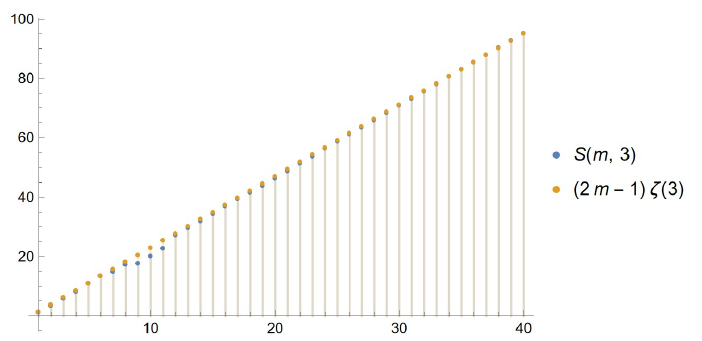}
		\caption{The full sum $S(m,3)$ follows the straight line $(2 \,m-1)\,\zeta(3)$ except for a slight deviation observed in the small interval $8\le m \le 12$. In that interval, a few blue dots are visible because they are not on the straight line. The cause of this deviation is explained in the next section. When the two functions completely overlap, the blue dots are hardly visible.}
\end{figure}
 \begin{figure}
	\centering
		\includegraphics[scale=0.60]{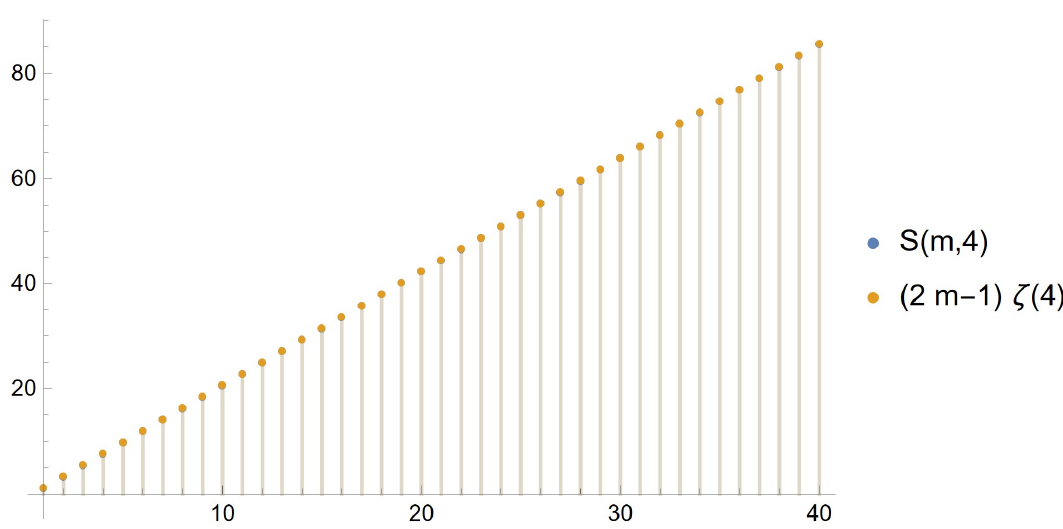}
		\caption{The full sum $S(m,4)$ and the straight line $(2 \,m-1)\,\zeta(4)$ overlap so exactly that they are indistinguishable on the plot i.e. the blue dots are covered almost completely by the orange/yellow dots.}
\end{figure}   
\begin{figure}
	\centering
		\includegraphics[scale=0.60]{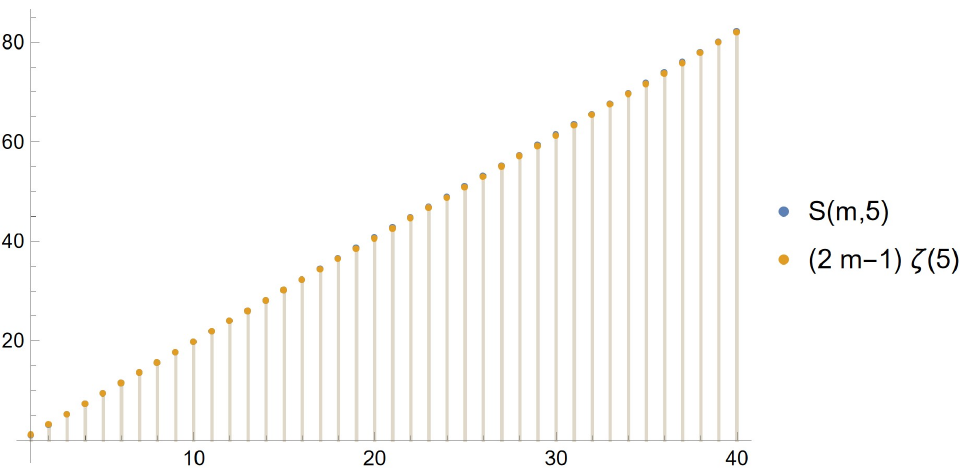}
		\caption{The full sum $S(m,5)$ and the straight line $(2 \,m-1)\,\zeta(5)$ overlap so exactly that they are indistinguishable on the plot i.e. again, the blue dots are covered almost completely by the orange/yellow dots.}
\end{figure}  

\section{The remainder $R(m,N)$}
In the previous section we split the full sum $S(m,N)$ into the dominant contribution $D(m,N)$ and the remainder $R(m,N)$. $D(m,N)$ could be evaluated analytically and was given by $(2 \,m-1)\,H_m^{(N)}$ (here $H_m^{(N)}=\sum_{x=1}^m 1/x^N$ approaches $\zeta(N)$ quickly after summing a few terms). In other words, $D(m,N)$ was essentially a straight line with slope $2\,\zeta(N)$. In Figures 1, 2 and 3 we plotted $S(m,N)$ and $D(m,N)$ on the same graph for $N=3$, $4$ and $5$ respectively. The two plots were both straight lines and overlapped almost exactly except for a small interval $8\le m \le 12$ in the $N=3$ case where $S(m,3)$ deviated slightly from a straight line. The near-equivalence of $S(m,N)$  and $D(m,N)$ implies that the remainder $R(m,N)$, given by the complicated sum \reff{Rmn}, must be negligible compared to $D(m,N)$ i.e. $|R(m,N)/D(m,N)|<<1$. There are therefore two possibilities: either the magnitude of $R(m,N)$ grows much slower than $2\,m$ (diverges but slowly) or converges to a finite value. Our numerical as well as analytical study below show that $R(m,3)$ likely diverges logarithmically whereas $R(m,N)$ converges to a finite value for $N\ge 4$. We also explain why $S(m,3)$ had a slight deviation from a straight line at $N=3$ in the small interval $8\le m \le 12\,$.

\subsection{Numerical plots of $R(m,N)$ and the special case of $N=3$} 
We begin by rewriting here the expression \reff{Rmn} for $R(m,N)$ for quick reference:
\beq
R(m,N)=\sum_{\substack{x,y,z=1\\z\ne x\,;\,z\ne y}}^m \dfrac{1}{x^N+y^N-z^N}\,.
\eeq{RmnAA}  
The first thing to note is that case $N=3$ is distinct from the other cases ($N>3$) because for $z\ne x$ and $z\ne y$ it is the only case where there exists positive integers $(x,y,z)$ where $f=x^N+y^N-z^N=\pm 1$. At these points, $1/f$ reaches the maximum magnitude possible of unity. There is an infinite number of points $(x,y,z)$ where $f=x^3+y^3-z^3=\pm 1$ since algebraic parametrizations yielding an infinite sequence of solutions have been known for a long time \cite{Mordell}. In a cubic region where $x$, $y$ and $z$ each run from $1$ to $m$, the number of such occurrences has been proven to have a lower bound of the order of $\mathcal{O}(m^{1/4})$ (see \cite{Mordell} for details and original references). We cannot deduce from this alone that $R(m,3)$ diverges since the sum involves positive and negative contributions. However, if we summed the absolute value of the terms, that is $|1/f|$ instead of $1/f$, the sum would clearly diverge since there would be an infinite number of $+1$ to sum. Therefore we know that the series at $N=3$ is not absolutely convergent. Moreover, the series would diverge roughly in the order of $\mathcal{O}(m^{\alpha})$ where $\alpha\approx 1/4<1$ so that the ratio with $D(m,N)$ given by approximately $m^\alpha/(2\, m-1)$, would still tend to zero as $m$ increases. So even in this worst case scenario where we sum the absolute value of the terms, the linearity of $S(m,3)$ would still persist. It then follows that the sum $R(m,3)$, where the terms $1/f$ can make a positive or negative contribution, would grow either more slowly or even possibly converge. The ratio with $D(m,N)$ would then be even smaller and the linearity of $S(m,3)$ would be even more robust. Whether $R(m,3)$ diverges or not is an interesting question but it has no bearing on the linearity of $S(m,3)$. As previously mentioned, our analytical study will suggest that $R(m,3)$ likely diverges logarithmically. 

We plot below $R(m,N)$ for $N=3$ and $N=4$ up to $m=140$ and for $N=5$ and $N=6$ up to $m=100$. We include a data table for $N=3$ and $N=4$. The numerical plots (figures 4-7) were generated with a regular laptop computer and the data tables with a powerful online computer simulation.  Of the four cases, the case $N=3$ fluctuates the most and also reaches the highest magnitude in the range $1\le m \le140$. The plot of $R(m,3)$ in Fig. 4 could not confirm whether it converged or not so a data table was generated at higher values of $m$ up to $3200$. However, even with the data table, it was not possible to determine whether $R(m,3)$ converged or not. The only thing one could conclude is that if there is a growth, it would be extremely slow.  The plot of $R(m,4)$ in Fig. 5 fluctuated considerably also.  We therefore generated a data table at higher $m$ values up to $3000$. This time, this showed that $R(m,4)$ does converge to a value close to $-0.029$. The plots of $R(m,5)$ and $R(m,6)$ in figures 6 and 7 respectively show they converge very quickly, with $N=5$ plateauing to $0.1673$ before reaching $m=40$ and $N=6$ plateauing or converging to $-0.00553$ before reaching $m=20$. $R(m,N)$ converges faster as $N$ increases. 

Note that there is a jump in the case $N=3$ at $m=9$. This occurs when $(x,y,z)$ is equal to $(6,8,9)$ (and $(8,6,9)\,$) where $6^3+8^3-9^3=-1$. This makes a large contribution of $-2$ to $R(m,3)$ so that it reaches $-2.78$ at $m=9$. At $m=12$, $R(m,3)$ reaches $\approx -0.5$ because when $(x,y,z)$ is equal to $(9,10,12)$ (and $(10,9,12)\,$) we obtain $9^3+10^3-12^3=+1$. This makes a large contribution of $+2$ to $R(m,3)$ so that it reverts back to a value close to $-0.5$ at $m=12$. We encounter $x^3+y^3-z^3=\pm 1$ next only at $m=103$ where $64^3+94^3-103^3=+1$. This 
shows up in Fig. 4 as a jump of $+2$ at $m=103$. The departure of $S(m,3)$ from a straight line occurs only in the region of small $m$ between $m=8$ and $m=12$. The jumps in $R(m,3)$ at larger $m$, such as at $m=103$, do not perturb much the linearity of $S(m,3)$ because at larger values of $m$, the value of the linear $D(m,3)$ term has grown very large compared to $R(m,3)$ that the ratio $|R(m,3)/D(m,3)|$ is insignificant. At $m=9$, the ratio is $0.14$ whereas at $m=103$, the ratio is $0.01$, a difference of a factor of $14$. The ratio of $0.14$ explains why the deviation of $S(m,3)$ from a  straight line occurs near $m=9$ but also why it is a slight and not drastic deviation.     
\begin{figure}
\centering
\begin{subfigure}{.5\textwidth}
  \centering
  \includegraphics[width=1.0\linewidth]{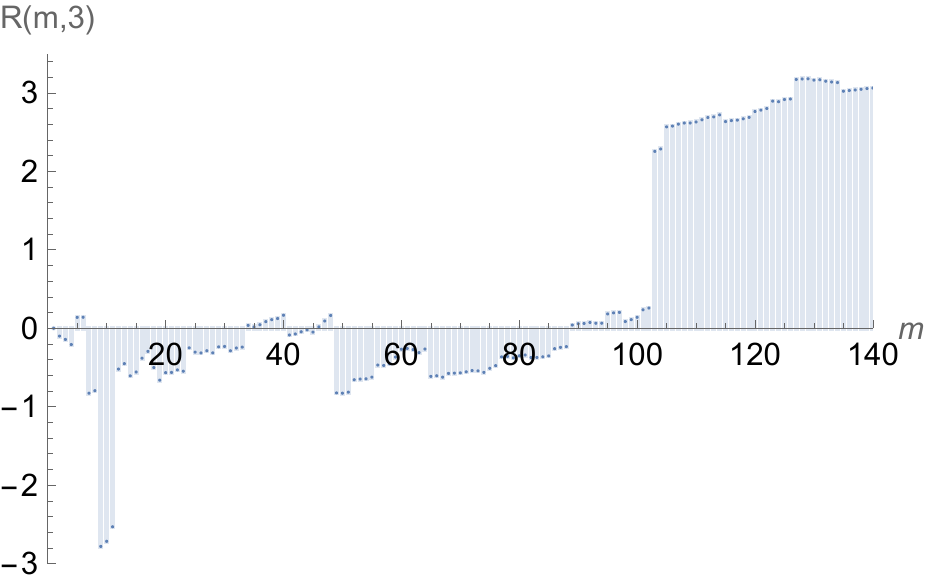}
  \caption{Plot of $R(m,3)$ vs. $m$ up to $m=140$. This has significant fluctuations. Note the jump at $m=9$ and it retreating back in magnitude at $m=12$. This corresponds to the contribution of Fermat near misses where $x^3+y^3-z^3=\pm 1$. The region near $m=9$ is where $S(m,3)$ had a slight deviation from a straight line in Fig. 1. The adjacent table of values goes up to $m=3200$ but one cannot deduce from it whether $R(m,3)$ converges or diverges.}
  \label{fig:sub1}
\end{subfigure}%
\begin{subfigure}{0.5\textwidth}
  \centering
  \includegraphics[width=0.7\linewidth]{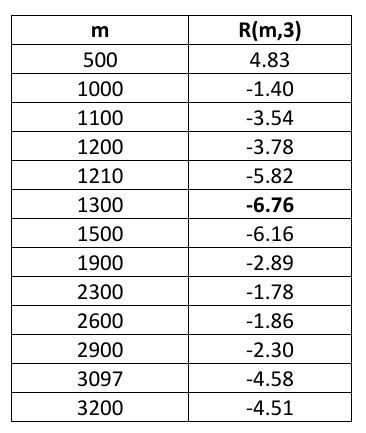}
  \caption{Table of values}
  \label{fig:sub2}
\end{subfigure}
\caption{$R(m,3)$}
\label{fig:test}
\end{figure}

\begin{figure}
\centering
\begin{subfigure}{.5\textwidth}
  \centering
  \includegraphics[width=1.0\linewidth]{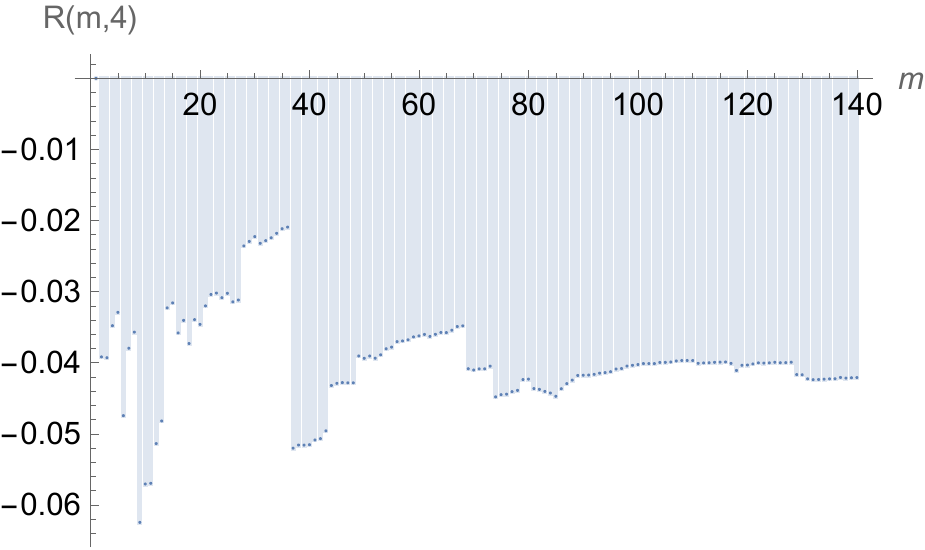}
  \caption{Plot of $R(m,4)$ vs. $m$ up to $m=140$. It shows significant fluctuations up to $m=140$ and one cannot tell if it converges or not until one goes to higher $m$ values. The adjacent table of values goes up to $m=3000$ and shows that it does eventually converge to $-0.029$.}
  \label{fig:subA1}
\end{subfigure}%
\begin{subfigure}{0.5\textwidth}
  \centering
  \includegraphics[width=0.7\linewidth]{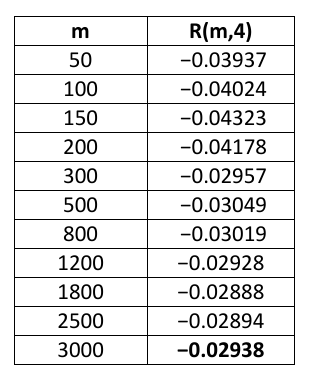}
  \caption{Table of values}
  \label{fig:subA2}
\end{subfigure}
\caption{$R(m,4)$}
\label{fig:Atest}
\end{figure}
 
\begin{figure}
	\centering
		\includegraphics[scale=0.60]{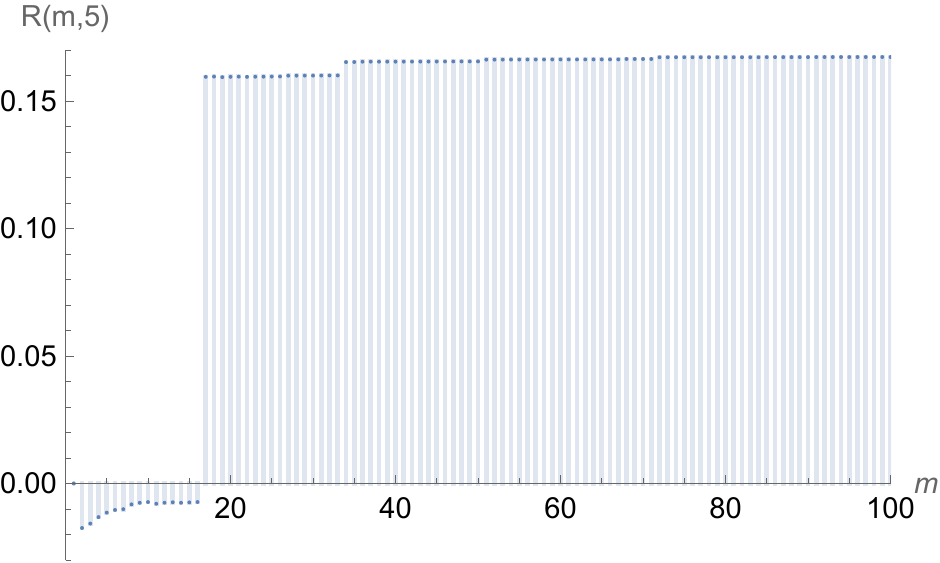}
		\caption{Plot of $R(m,5)$ vs. $m$  up to $m=100$. It converges quickly (before reaching $m=40$) to the value of $0.1673$. So the case $N=5$ converges significantly faster than the case $N=4$. This trend continues: as $N$ increases, the convergence is faster.}
\end{figure}  
\begin{figure}
	\centering
		\includegraphics[scale=0.60]{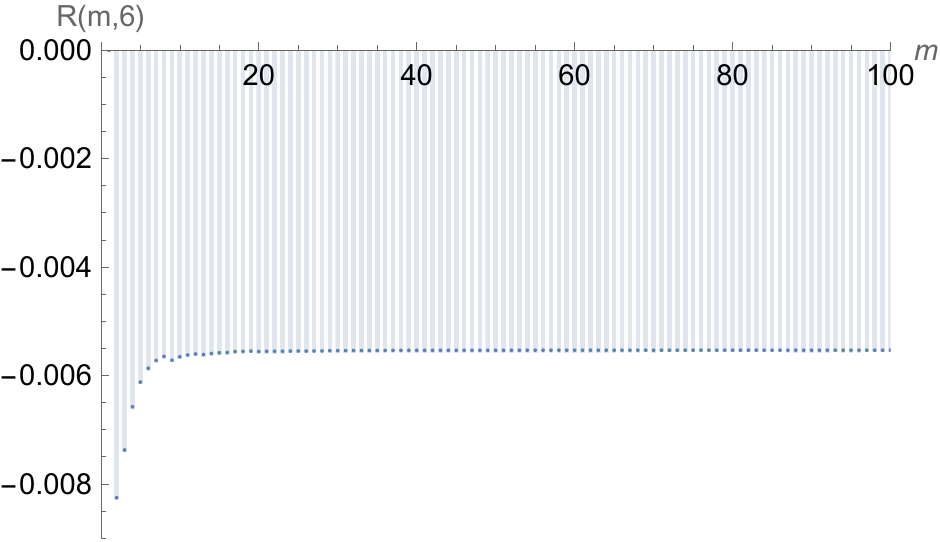}
		\caption{Plot of $R(m,6)$ vs. $m$  up to $m=100$. It converges very quickly (before reaching $m=20$) to the value of $-0.00553$. So the case $N=6$ converges faster than the case $N=5$ and its magnitude is smaller.}
\end{figure} 

\clearpage
\subsection{Analytical study of convergence: comparing the sum $R(m,N)$ to the Cauchy principal value of its integral} 
The previous numerical analysis suggests that for $N\ge4$, $R(m,N)$ converges. However, the case of $R(m,3)$ was not as clear (at least from the numerical analysis we carried out). In this section, by comparing the sum $R(m,N)$ to the Cauchy principal value of its integral, we find that the remainder $R(m,3)$ is expected to diverge logarithmically and that $R(m,N)$ for $N\ge 4$ is expected to converge. Since we are not applying here the traditional integral test (our integrand does not meet its criteria), this cannot be viewed for now as a rigorous proof. Why it is reasonable is discussed below.

The remainder $R(m,N)$ is defined by \reff{Rmn} which we rewrite here for quick reference:
\beq
R(m,N)=\sum_{\substack{x,y,z=1\\z\ne x\,;\,z\ne y}}^m \dfrac{1}{x^N+y^N-z^N}\,.
\eeq{RN} 
The original sum $S(m,N)$ does not have the constraint $z\ne x$ and $z\ne y$. The contributions at $z=x$ and $z=y$ of the original sum are large and important as this is what yields the linearity in $m$ (which diverges as $m\to \infty$). The non-trivial contributions of $z=x$ and $z=y$ to the original triple sum does not have an integral counterpart. If $S(m,N)$ were turned into a triple integral, the contribution to the integral from say $z=x$ would be zero as it would correspond to a plane with no thickness and hence zero volume. However, the remainder $R(m,N)$ specifically excludes $z=x$ and $z=y$. Nonetheless, one must address two more potential obstacles which are of a number-theoretic nature before one can establish an integral or continuum version for $R(m,N)$. One is the fact that the denominator in $R(m,N)$ is never zero due to Fermat's Last Theorem whereas in an integral the denominator would reach zero. This is not an issue if we take the Cauchy principal value of the integral (where as usual the point that yields zero is bypassed by integrating up to a distance of $\epsilon$ away from it on both sides so that a cancellation occurs and the limit as $\epsilon\to 0$ yields a well-defined result). The second obstacle are the Fermat near misses. This is particularly important in the case of $N=3$ where there exists infinitely many solutions to $x^3+y^3-z^3=\pm 1$ for $z\ne x$ and $z \ne y$ (these are called non-trivial solutions in contrast to the trivial ones when $z=x$ or $z=y$). There is no integral analog to these Fermat near misses. However, if we assume that the $+1$ and $-1$ are equally distributed up to large (infinite) $m$, the sum of the positive and negative contributions would basically cancel out or be negligible. There does not appear to be a theorem that has been proved that states that the $+1$ and $-1$ are equally distributed when $x$, $y$ and $z$ are positive integers with the constraint $z\ne x$ and $z\ne y$. However, numerical simulations seem to suggest that: at $m=3000$ there is zero difference, at $m=5000$ there is an imbalance of $-2$ ($-1$ occurs two more times) and at $m=10,000$ there is again a difference of zero. So the difference between the occurrences of $+1$ and $-1$ appears to never grow, is extremely low or non-existent. Even though $+1$ or $-1$ Fermat near misses are absent in the case of $N\ge 4$ (for $z\ne x$ and $z\ne y$) we expect that any accumulation due to any possible degeneracies will first of all be more sparse and moreover have positive and negative contributions that also roughly balance out. 

We are therefore ready to compare the Cauchy principal value (PV) of the integral to the sum $R(m,N)$. The PV integral analog to $R(m,N)$ with $m\to \infty$ is labeled $I_{PV}$ and is given by
\beq
I_{PV} =\int_{x=1}^{\infty} \int_{y=1}^{\infty} PV \int_{z=1}^{\infty} \dfrac{1}{x^N+y^N-z^N}\,dx\,dy\,dz\,.
\eeq{IPN}            
The Cauchy principal value of the integral over $z$ yields the following function of $x$ and $y$:
\begin{align}
S(x,y)&=\frac{\pi  \cot \left(\frac{\pi }{N}\right)}{N \left(x^N+y^N\right)^{\frac{N-1}{N}}}-\frac{B_{(x^N+y^N)^{-1}}\left(\frac{1}{N},0\right)}{N \left(x^N+y^N\right)^{\frac{N-1}{N}}}
\label{PVN}
\end{align}
where $B_z(a,b)$ is the incomplete Beta function defined as 
\beq
B_z(a,b)=\int_0^z t^{a-1}\,(1-t)^{b-1}\,dt\,.
\eeq{Bz}
We now need to carry out the integral over $x$ and $y$ of $S(x,y)$. It is best to go over to polar coordinates as the convergence or divergence of the integral is determined by the behavior of the integrand in the large $r$ limit. We define $x= r \,\cos(\theta)^{2/N}$ and $y=r\,\sin(\theta)^{2/N}$. Since both $x$ and $y$ are positive, the angle $\theta$ ranges between $0$ and $\pi/2$ but does not include $0$ or $\pi/2$ since $x\ge 1$ and $y\ge 1$ (the origin $(x,y)=(1,1)$ is located at $\theta=\pi/4$). We now have that 
$x^N+y^N =r^N$. The Jacobian for the transformation is $J=f(\theta)\,r$ where $f(\theta)$ is a function of $\theta$ that can easily be determined but plays no role in determining convergence or divergence. The function $S(x,y)$ has no dependence on $\theta$  and is given by 
\beq
S(r)=\frac{\pi \, \cot \left(\frac{\pi }{N}\right)}{N \,r^{N-1}}-\frac{B_{r^{-N}}\left(\frac{1}{N},0\right)}{N \,r^{N-1}}\,.
\eeq{Sr}     
Therefore the integrand $U$ is given by $r\,S(r)\,f(\theta)$ which is given by
\beq
U=\frac{f(\theta)\,\pi\, \cot \left(\frac{\pi }{N}\right)}{N \,r^{N-2}}-\frac{f(\theta)\,B_{r^{-N}}\left(\frac{1}{N},0\right)}{N \,r^{N-2}}\,.
\eeq{II}
The incomplete Beta function about large (infinite) $r$ behaves as $B_{r^{-N}}\left(\frac{1}{N},0\right)\approx \frac{N}{r}$. Therefore the integrand $U$ in the large $r$ limit behaves as
\beq
U_{r\to \infty} \approx \frac{f(\theta)\,\pi\, \cot \left(\frac{\pi }{N}\right)}{N \,r^{N-2}}-\frac{f(\theta)}{r^{N-1}}\,. 
\eeq{Irr}
There are two terms: the first term behaves as $1/r^{N-2}$ and the second term as $1/r^{N-1}$. For $N=3$, the first term falls off as $1/r$ and leads to a logarithmic divergence when integrated over r (whereas the second term yields $1/r^2$ and is finite when integrated). Therefore, $R(m,3)$ is expected to diverge logarithmically, that is as $\log(r)$ which is roughly equivalent to saying that $R(m,3)$ diverges as $\log(m)$. At $N=4$, the first term falls off as $1/r^2$ and the second term falls off as $1/r^3$. When integrated over $r$, $R(m,4)$ therefore falls off as $1/r$ and hence converges in agreement with our numerical analysis. As $N$ increases further, $R(m,N)$ falls off faster at a rate of $1/r^{N-3}$ so that all $N \ge 4$ cases converge. We therefore expect a faster convergence for $N=5$ than for $N=4$ and a faster convergence for $N=6$ than for $N=5$. That is precisely what we observed numerically in figures $5$, $6$ and $7$ representing $N=4$, $N=5$ and $N=6$ respectively. 

We have applied a test involving the Cauchy principal value (PV) of an integral to the sum $R(m,N)$. Note that such a test was applied to a sum $R(m,N)$ which had constraints, namely $z\ne x$ and $z\ne y$. Without such constraints, the test would clearly not be valid. The only thing that can change our conclusions of this section from being plausible to a more rigorous proof is a theorem on a Cauchy principal value integral test that can be used when the integrand is not strictly positive, encounters singularities that can be dealt with by PV, and is not strictly decreasing. Such a theorem would apply to sums obeying certain criteria (e.g. $S(m,N)$ would not pass those criteria).

\section{Conclusion}
In this paper we investigated the sum of the reciprocal of the polynomial $f=x^N+y^N-z^N$which appears in Fermat's Last Theorem. Here $N\ge 3$ is an integer and $x$, $y$ and 
$z$ are positive integers. The interesting thing here is that the polynomial can be positive or negative but is never equal to zero due to Fermat's Last Theorem. We therefore have a well-defined sum that includes positive and negative terms. The value of performing such a sum is that its behavior is remarkably simple. The finite sum $S(m,N)$ of the reciprocal $1/f$ over $x$, $y$ and $z$ that run from $1$ to $m$ inclusively was split into two parts: the dominant contribution $D(m,N)$ and a remainder $R(m,N)$. The dominant contribution could be determined analytically and was given  by $D(m,N)=(2\,m-1)\,H_m^{(N)}$ where $H_m^{(N)}=\sum_{x=1}^m 1/x^N$. $H_m^{(N)}$ approaches quickly the Riemann zeta function $\zeta(N)$ after just a few terms. Therefore $D(m,N)$ was linear in $m$ with slope $2\,\zeta(N)$. The remainder $R(m,N)$ was a complicated sum but it was completely negligible compared to $D(m,N)$. Therefore $S(m,N)$ was almost identical to $D(m,N)$. In figures 1, 2 and 3 we plotted $S(m,N)$ vs. $m$  for the cases of $N=3,4$ and $5$ respectively and all of them were straight lines (with a  small caveat for $N=3$). The dependence of the slope on $N$ is weak since $\zeta(N)$ does not change significantly with $N$ and has a tight range between $1.202$ (at $N=3$) and unity (at large $N$). Not only is the original sum $S(m,N)$ a straight line but the slope hardly varies with $N$. Instead of making things more complicated, the negative sign in front of $z^N$ has led to a great simplification. 

The case $N=3$ is the only case where a slight deviation from a straight line was observed for $S(m,N)$; this took place over a small interval near $m=9$. We were able to show that this was due to Fermat near misses at small $m$ where $x^3+y^3-z^3=\pm 1$. A jump could be observed in the remainder $R(m,3)$ in the small interval near $m=9$ where the deviation occurred.  

We made numerical plots of the remainder $R(m,N)$ vs. $m$ for $N=3,4,5$ and $6$ and these appear in figures $4$, $5$, $6$ and $7$. We also generated a table of values for $R(m,3)$ and $R(m,4)$ that went up to larger $m$ values (around $3000$) since the plots up to $m=140$ for these two cases did not converge or plateau to a particular value in contrast to $N=5$ and $N=6$ that converged to a value quickly. The case $N=3$ had the most fluctuations and the table of values that went up to $m=3200$ was not sufficient to determine if it diverged or not. By making an analytical study of the Cauchy principal value of its integral version, we found it likely diverges logarithmically, that is as $\log(m)$. $R(m,4)$ fluctuated a lot also but in contrast to $R(m,3)$, the table of values showed that it ultimately converged, albeit very slowly. Our analytical study of the Cauchy principal value of the integral version of $R(m,N)$ showed that $R(m,N)$ converged for $N\ge 4$ and moreover, converged faster as $N$ increased, in agreement with our numerical analysis.  

We determined whether $R(m,N)$ converges by comparing it to the Cauchy principal value of its integral. We explained why such a comparison was reasonable but a theorem that establishes a Cauchy principal value (PV) integral test would nonetheless be welcome. This could be applied to test the convergence of a sum where the integrand of its integral version is not strictly positive, has singularities that can be resolved (by evaluating the PV) and is not strictly decreasing in each variable. There is a separate interesting question that is worth exploring: if the absolute value of the terms in the sum $R(m,3)$ were used instead i.e. the sum of $1/|(x^3+y^3-z^3)|$ with the same constraint $z\ne x$ and $z\ne y$, would this modify the logarithmic divergence of $R(m,3)$ and if so, how? The Fermat near misses of $+1$ and $-1$ would now each contribute $+1$ to the sum. This could make a substantial contribution to the sum depending on how dense the near misses are and this could result in a clear departure from a logarithmic divergence. Similarly, we found that $R(m,N)$ for $N \ge 4$ converges but we did not determine whether they converged absolutely. A preliminary numerical analysis suggests they do. It would therefore be worthwhile to find an analytical proof of this.

\end{document}